\documentclass[,a4paper,12pt]{article}
\usepackage{amsfonts, amssymb, amscd, amsmath}
\input{xy}
\usepackage[english]{babel}
\usepackage[all]{xy}

\newtheorem{theo}{Theorem}[section]
\newtheorem{coro}[theo]{Corollary}
\newtheorem{lemm}[theo]{Lemma}
\newtheorem{prop}[theo]{Proposition}
\newtheorem{defi}[theo]{Definition}
\newtheorem{rema}[theo]{Remark}
\newtheorem{exam}[theo]{Example}

\newenvironment{proo}{\noindent \textbf{{Proof.}} \sf}

\def\qed{\hfill $\diamond$ \bigskip}

\newcommand{\Bdf}{\begin{defi}}
\newcommand{\Edf}{\end{defi}}
\newcommand{\Bte}{\begin{theo}}
\newcommand{\Ete}{\end{theo}}
\newcommand{\Bpo}{\begin{prop}}
\newcommand{\Epo}{\end{prop}}
\newcommand{\Bcr}{\begin{coro}}
\newcommand{\Ecr}{\end{coro}}
\newcommand{\Blm}{\begin{lemm}}
\newcommand{\Elm}{\end{lemm}}
\newcommand{\Bex}{\begin{exam}}
\newcommand{\Eex}{\end{exam}}
\newcommand{\Bdm}{\begin{proo}}
\newcommand{\Edm}{\end{proo}}
\newcommand{\Brm}{\begin{rema}}
\newcommand{\Erm}{\end{rema}}

\newcommand{\li}[2]{{\vphantom{#2}}_{#1}#2}
\def\im{{\mathsf{Im}}}

\def\B{{\mathcal B}}
\def\C{{\mathcal C}}

\def\F{{\mathcal F}}

\def\K{{\mathcal K}}
\def\L{{\mathcal L}}

\def\M{{\mathcal M}}
\def\N{{\mathcal N}}

\def\X{{\mathcal X}}

\begin{document}
\sf
\title{\sf Invariants of a Free Linear Category and Representation Type}
\author{Claude Cibils and Eduardo N. Marcos
\thanks{\sf The authors were partially supported by the project Mathamsud GR2HOPF. The second named author  was partially supported by Conselho Nacional de Desenvolvimento Cient\'{\i}fico
Tecnol\'{o}gico (CNPq, Brazil) and by the Projeto tem\'atico FAPESP 2014/09310-5 }}
\date{}
\maketitle

\begin{abstract}
We consider a homogeneous action of a finite group on a free linear category over a field in order to prove that the subcategory of invariants is still free. Moreover we show that the representation type is preserved when considering invariants.
\end{abstract}


\section{\sf Introduction}
The first purpose of this article is to prove that the category of invariants of a homogeneous action of a finite group on a free linear category is again a free linear category. V.K. Kharchenko \cite{kha} and D.R.  Lane \cite{lan} proved  that the algebra of invariants of a finite group acting homogeneously on a free algebra is a free algebra.

In this paper $k$ is a field of any characteristic. A  $k$-category is a small category enriched over $k$-vector spaces, see for instance \cite{ke}. In other words objects are a set, morphisms are vector spaces and composition is $k$-bilinear. At each object the endomorphisms forms a $k$-algebra.

A free $k$-category is given  by a set of objects and a set of  "oriented" vector spaces, that is to each couple of objects there is a given  vector space which can be zero. The  \textit{track-quiver} of this data is the oriented graph which records the non-zero vector spaces. In the first Section we recall the precise definitions and the construction of the associated free $k$-category. We also recall that if the number of objects of a $k$-category is finite, there is a canonical $k$-algebra associated to it. This algebra is hereditary if the category is free.

The problem of describing invariants of an action of a group goes back to D. Hilbert and E. Noether. For commutative free algebras (i.e. polynomial algebras) over a field of characteristic zero, a homogeneous action of a finite subgroup of $GL(d,k)$ provides again a commutative free algebra of invariants if and only if the group is generated by pseudo-reflections. This has been proved by G.C. Shephard -- J.A. Todd \cite{shetod}, and by C. Chevalley and J.P. Serre \cite{che,ser}. For a detailed account see \cite{sta,ben}.

The action of a group on a free $k$-category is called homogeneous if it is given by a linear action on the generators,  which is extended in the unique possible manner to an action on the free $k$-category by invertible endofunctors. Note that we consider actions on a free $k$-category which are trivial  on the set of objects. Indeed, we restrict to this case since otherwise the invariants are the invariants of the full subcategory determined by the invariant objects. Observe that this context  is at the opposite of a Galois action of a group on a $k$-category where the action on objects has to be free, hence with no invariant objects, see for instance \cite{cibmar}.

Another aim of this paper is to prove a result about invariants of a tensor product of $kG$-modules which appears to be unknown. This result is a crucial tool for proving our main theorem, namely that invariants of a homogeneous action on a free $k$-category form again a free $k$-category. It can be outlined as follows. Let $M_n,\dots,M_1$ be a sequence of $kG$ modules and let $\left(M_n\otimes\dots\otimes M_1\right)^G$ be the vector space of  invariants of the action of $G$ on the tensor product. Some of these invariants  "are strictly from below" in the sense that they are sums of tensors of invariants of a  partition by strict sub-strings of  $M_n\otimes\dots\otimes M_1$. We call those invariants "composites", they form a canonical sub-vector space of the space of invariants. An arbitrarily chosen vector space complement is always called a space of irreducible invariants. For any fixed choice of irreducible spaces of invariants in the sub-strings of the tensor product, we prove  the uniqueness of the decomposition of an invariant as a sum of tensors of irreducible invariants. We infer this result from the Kharchenko--Lane Theorem.

The last section is also a contribution to the classical study of the relationship of an algebra with its algebra of invariants, see for instance \cite{reri}. We prove that the representation type is conserved, more precisely if a free $k$-linear category with a homogeneous action of a finite group is of finite or tame representation type, then its  invariant $k$-category is respectively of finite or tame representation type. Note that we adopt the convention that a $k$-category of finite representation type is also of tame representation type.  In order to prove this result, we set up ad-hoc cleaving techniques whose basics goes back to D.G. Higman.

At least two questions arise from the present work in relation with previous results for free algebras.

Firstly V.O. Ferreira, L.S.I. Murakami and A. Paques in \cite{fermurpaq} have provided a generalization of the work of V.K. Kharchenko by considering a homogeneous action of a Hopf  $k$-algebra over a free $k$-algebra. A generalization of our work  should also hold for a Hopf algebra acting homogeneously on a  free $k$-category.

Secondly note that W. Dicks -- E. Formanek \cite{dicfor} have proved that for a free $k$-algebra, unless the group is cyclic and the number of free generators of the original algebra is finite, the number of free generators of the invariants is infinite. A similar result in the setting of free $k$-categories should hold.

Acknowledgement: We thank the referee for suggestions which have improved the paper, for detecting imprecise or erroneous statements as well as for an attentive reading.

\section{\sf Free linear categories}
We recall first the definition of a free linear category ${\L}_k(V)$ over a field $k$ which satisfies the universal property stated in Proposition \ref{up}..  For this purpose, let   $\L_0$ be a set and let $V=\{_y\!V_x\}_{y, x \in \L_0}$ be a family of $k$-vector spaces. Let $u=(u_n, \cdots, u_0)$ be a sequence of elements in $\L_0$ and consider the  vector space
$$W(u) =\li{u_n}{\!V}_{u_{n-1}}\otimes \cdots \otimes\li{u_2}{\!V}_{v_1}    \otimes  {}_{u{_1}}\!V_{u_{0}}.$$
 For a singleton sequence $u=(u_0)$ we set $W(u_0) =k$.

 Let ${}_y{S}_x$ be the set of sequences $(y=u_n, u_{n-1}, \cdots, u_0= x)$ and let
 $$ \li{y}{W}_x = {\bigoplus_{u \in {\li{y}{S}_x}}} W(u).$$

 Notice that the concatenation of sequences of elements of $\L_0$
  $$v=(z=v_m, \dots, v_0=y)\in \li{z}{S}_y \mbox{\ and\ } u=(y=u_n,\dots , u_0 =x)\in \li{y}{S}_x$$  is the sequence
 $$vu =(z=v_m, \dots, v_0=y=u_n, \dots ,u_0=x).$$

\begin{defi}\label{definition of free} Let $\L_0$ be a set, and let $V$ be a family as above. The \emph{free $k$-category $\L_k(V)$} has set of objects $\L_0$. For  $x, y\in \L_0$, the vector space of morphisms $_y\left[\L_k(V)\right]_x$ from $x$ to $y$ is ${}_y W_x$.
The composition is given by the direct sum of the tensoring linear maps $$W(v)\otimes W(u) \to W(vu).$$.
 \end{defi}

We observe that the preceding definition has the advantage of being intrinsic, since it do not rely on the choice of  bases of the vector spaces of the family. We recall the universal property of $\L_k(V)$:

\Bpo \label{up} Let $\B$ be a $k$-category with set of objects $\B_0$. A linear functor $F:\L \to \B$ is uniquely determined by the following:
\begin{itemize}
\item[-] a map
$F_0: \L_0 \to \B_0$,
\item[-]
 a family of linear maps $\{{}_yF_x:\, _y\!V_x \ \to\  {}_{{F_0}(y)}\B_{{F_0}(x)}\}_{x,y\in \B_0}.$
\end{itemize}
\Epo

 To the data required for defining a free $k$-category we associate a \textit{track-quiver}, that is an oriented graph whose vertices are $\L_0$ and with an arrow from $x$ to $y$ if and only if $\li{y}{V}_x \neq 0$.
 We also associate a \textit{bases-quiver} as follows: firstly choose a basis of each vector space in the family. The vertices are still $\L_0$ while the set of arrows from $x$ to $y$ is the chosen basis of  $\li{y}{V}_x $. Those arrows are said to have \textit{source} $x$ and \textit{target} $y$.

\begin{rema}\label{free}
Let $\C$ be a  small category (i.e. without additional $k$-structure). By definition its linearization $k\C$ is a $k$-category with the same set of objects, while morphisms are the vector spaces with bases the morphisms of the category. Composition is obtained by extending bilinearly the composition of $C$.

  Associated to a quiver $Q$ there is a free category $\F_Q$ (without additional structure) as follows, see for instance \cite[p. 48]{maclane}: the objects are the vertices, morphisms are sequences of concatenated arrows (called paths), including the trivial one at each vertex. The source (resp. target) of a non-trivial path is the source (resp. target) of its first (resp. last) arrow. The source and target of a trivial path at a vertex are that vertex. Composition is concatenation of paths and the trivial paths are the identities.

 Consider a free $k$-category $\L_k(V)$ with bases-quiver $Q$. The linearization  $k\F_Q$ is isomorphic to $\L_k(V)$.
\end{rema}
 We recall now  the definition of the $k$-algebra corresponding to a $k$-category with a finite number of objects.

 \begin{defi} Let $B$ be a $k$-category with set of objects $B_0$. The $k$-algebra $\oplus\B$ is the vector space $\displaystyle \oplus_{x, y\in\B_0}\
{_xB_y}$ equipped with the product extended linearly from the following: if $g \in {}_zB_{y'}$ and $f\in {}_yB_x$, then $gf= 0$ if $y\neq y'$ while if $y=y'$ the product is the composition $gf$ of the morphisms in the category.
 \end{defi}
 Note that the algebra above has a unit if and only if the set $B_0$ is finite, in that case we say that $\B$ is \textit{object-finite}. The following result is clear.

 \Blm
  Let $\L_0 =\{x\}$ be a singleton, let the family $V$ be constituted by a single vector space $\{{} _x V_x\}$ and let $\L_k(V)$ be the corresponding free $k$-category. Then $\oplus{\L_k(V)} =
T_k({} _x\!V_x)$ where the latter is the tensor $k$-algebra
$$k\oplus {} _x V_x\oplus \left({} _x V_x\otimes_k{} _x V_x\right)\oplus\cdots \oplus\left({} _x V_x\right)^{\otimes_k^n}\oplus\cdots $$
\Elm
 We recall next the following  well known result, see for instance \cite{cib}.

 \Bpo
  Let $\L_k(V)$ be a free object-finite $k$-category. Then $\oplus{\L_k(V)}$ is a hereditary algebra.
 \Epo
\Bdm
Let $E = \displaystyle{\times_{s\in \L_0}} W(s)$ be the semisimple subalgebra of $\Lambda = \oplus{\L_k(V)}$ defined by the objets of $\L$.
Let
$$M = \displaystyle{\bigoplus_{x,y \in \L_0}}\  _y V_x.$$
The following is an exact sequence of $\Lambda$ bimodules
$$ 0\longrightarrow \Lambda \otimes_E M \otimes_E \Lambda \stackrel{\alpha}{\longrightarrow} \Lambda
 \otimes_E \Lambda \stackrel{m}{\longrightarrow} \Lambda \longrightarrow 0$$
where $m$ is the multiplication given by $\Lambda$ and $\alpha(\lambda \otimes m\otimes \mu) = \lambda m\otimes \mu - \lambda \otimes m\mu.$
Note that $\Lambda \otimes_E M \otimes_E \Lambda$ and $\Lambda \otimes_E \Lambda$ are projective $\Lambda$ bimodules. So $\Lambda$ has projective dimension smaller or equal 1 as a
bimodule, that is $\Lambda$ is hereditary.
\Edm

\Brm Let $T_E(M) = E\oplus M \oplus (M\otimes_E M)\oplus\  \cdots \oplus (M^{\otimes^n_E}) \oplus \cdots$ be the tensor algebra of
the bimodule $M$ over $E$. Then $\oplus\L_k(V) = T_E(M)$. Moreover if we choose a basis of each $_yV_x$ and we consider $Q$ the corresponding
bases-quiver,  $T_E(M)$ is isomorphic to the path algebra $kQ$.
\Erm

We recall a well known fact, namely a sub-category of a free $k$-category is not free in general.
\begin{exam}\label{notfreesubcat}
Let $\L$ be the free $k$-category whose bases-quiver is
$$u_0\longrightarrow u_1 \longrightarrow u_2 \longrightarrow u_3$$
Note that the generating vector spaces are of dimension $1$. Consider the subcategory $\B$ given by
 $\B_0 = \L_0$ while  $_{u_2}\B_{u_1}=0$ and ${}_y\B_x = {}_y\L_x$ for all other pairs of objects $x$ and $y$. Then
$\B$ is the incidence $k$-category of the partially ordered set $$\{u_0,u_1,u_2,u_3 \ | \ u_0 < u_1, u_0< u_2,  u_1 < u_3, u_2<u_3\}.$$ It is well known and easy to prove that the algebra of this category have global dimension 2. Therefore the previous Proposition shows that $\B$ is not free.
\end{exam}
\section{\sf {Unique decomposition of invariants by irreducible}}
Let $G$ be a finite group, let $kG$ be the group algebra and let $V$ be a $kG$-module.
Since $kG$ is a Hopf algebra the tensor product over $k$ of two $kG$ modules $V$ and $V'$  is a $kG$-module through
the comultiplication \textit{i.e.} the diagonal action. More precisely for $s\in G, v\in V$ and $v'\in V'$
$$s(v\otimes v') = sv\otimes sv'.$$
We infer the action of $G$ on $T_k(V)$  by automorphisms of the algebra. The subalgebra of invariants (or the fixed points algebra)  is denoted $T_k(V)^G$.
\Bte (Kharchenko--Lane, see \cite{kha,lan})

Let $G$ be a finite group and let $V$ be a $kG$-module which is a finite dimensional $k$-vector space.
Then $T_k(V)^G$ is again a tensor algebra and there exists a homogeneous sub-vector space $U\subset T_k(V)^G$ such that $T_k(U)$ is isomorphic to
$T_k(V)^G$.
\Ete

\Brm The homogeneity of the free generators of $T_k(V)^G$ follows from the proof of V.K. Kharchenko.
\Erm

\begin{defi}\label{sg}
A free $k-$category $\L$ over a set $\L_0$ and over a family of $k$-vector spaces $V=\{{}_yV_x\}_{x,y \in \L_0}$ is called \emph{Schurian-generated} if each vector space of the family $V$ is one-dimensional or zero.
\end{defi}

Observe that a free $k$-category $\L_k(V)$ is Schurian-generated if and only if its bases-quiver has no multiple arrows, that is if the track-quiver and the basis-quiver are equal.

Next we will prove an analogous of the Kharchenko-Lane Theorem for free Schurian-generated $k$-categories. Later on we will prove a general theorem using other methods, however this result is useful as a first approach. Moreover Schurian-generated categories will be needed at the last Section.

\begin{prop}\label{sgfree}
Let $\L_k(V)$ be a free Schurian-generated $k$-category. Let $G$ be a finite group acting homogeneously on $\L_k(V)$. Then the invariant subcategory $\L_k(V)^G$ is again free.
\end{prop}

\begin{proo}
Since the bases-quiver has no multiple arrows,  to each arrow $a$ there is an associated character $\chi_a : G\to k^{\cdot}$. Let $w$ be a path of positive length in the bases-quiver and let $\chi_w$ be the product of the characters corresponding to the arrows of $w$. The path $w$ is an invariant if and only if $\chi_w = 1$.  The invariant paths from $x$ to $y$ form bases of ${}_y[\L_k(V)^G]_x$. An invariant path of positive length is called irreducible if it is not the concatenation of two invariant paths of positive length. Clearly any invariant path is a unique concatenation of irreducible invariant paths. Then $\L_k(V)^G$ is free over the set $\L_0$ and over the family of vector spaces with bases the invariant irreducible paths. Observe that its bases-quiver has an arrow corresponding to each irreducible invariant path, sharing the same source and target vertices.
\qed

\begin{exam}\label{crown}
Let $Q_n$ be an  oriented crown, that is a quiver with set of vertices being the set of elements of the cyclic group $C_n=<t \mid t^n=1>$ and an arrow from $t^i$ to $t^{i+1}$ for every $i$. Note that the corresponding free $k$-category is connected and that all the vector spaces of morphisms are infinite dimensional. Let $q$ be a primitive $n$-th root of unity and let the generator $t$ of $C_n$ act on each arrow of $Q_n$ by multiplication by $q$. The invariant irreducible paths are the paths of length $n$. The bases-quiver of the invariant $k$-category is a union of $n$ loops. In other words the invariant sub-algebra is a product of n-copies of the polynomial algebra in one variable. We observe that the invariant $k$-category is not connected.
\end{exam}
\end{proo}

From the Kharchenko-Lane Theorem we will infer a result about invariants in a tensor product of $KG$ modules which will be an important tool
in the next section.

Let $M_n, \dots, M_1$ be a sequence of $kG$-modules. Our next purpose is to distinguish  in
$(M_n\otimes \cdots \otimes M_1)^G$  between  ``composite invariants", namely invariant elements which comes strictly from below,  and  ``irreducible invariants".

For each natural number $i$ between 1 and $n-1$ consider the canonical map:
$$\varphi_i: (M_n\otimes \cdots \otimes M_{i+1})^G \otimes (M_i\otimes\cdots M_1)^G \longrightarrow (M_n\otimes \cdots \otimes M_1)^G.$$

\begin{defi} Let $\varphi = \displaystyle{\sum_{i=1}^{n-1}} \varphi_i$. We will call \emph{space of composite invariants}   the
image of $\varphi$  which we will  denote by $ \left[(M_n\otimes \cdots \otimes M_1)^{G}\right]^2 .$

\end{defi}

There is not a canonical complement to the subspace of composite invariants in the space of invariants. An arbitrarily chosen complement will be always called \textit{a space of irreducible invariants}.

Let $j,i$ be integers such that  $n\geq j\geq i\geq 1$. Then $M_j\otimes \cdots \otimes M_i$ will be called the $[j,i]-$\textit{string} and  $(M_j\otimes \cdots \otimes M_i)^G$ the $[j,i]-$\textit{invariant string}. We denote  $(M_j\otimes \cdots \otimes M_i)^{G}_{\mathsf{irr}}$ a chosen vector space which complements  $\left[(M_j\otimes \cdots \otimes M_i)^{G}\right]^2$ in the  $[j,i]-${invariant string}. Note that $(M_i)^G_{{\mathsf{irr}}} = M_i^G$
since $\left[M_i^{G}\right]^2 = 0$.

Let $p=(n_l, \dots, n_1)$ be a \textit{non-ordered $l$-partition} of $n$, that is a  sequence of $l$ positive integers such that $n=n_l+\cdots +n_1$.
 Let $(M_n \otimes\cdots \otimes M_1)^G_{p, \mathsf{irr}}$ be the following vector space
$$(M_n\otimes \cdots \otimes M_{n_1 + \dots n_{l-1} + 1})^{G}_{\mathsf{irr}} \otimes\cdots\otimes (M_{n_1+n_2}\otimes\cdots\otimes M_{n_1+1})^G_{\mathsf{irr}}\otimes (M_{n_1}\otimes \cdots \otimes M_{1})^{G}_{\mathsf{irr}}$$
and consider the map $$\psi_p : (M_n\otimes\cdots M_1)^G_{p, \mathsf{irr}} \to (M_n\otimes \cdots \otimes M_1)^G .$$ Let $\psi_{M_n,\dots,M_1}$ be the sum of the $\psi_p$ along the set all non-ordered partitions of $n$, namely
$$\psi_{M_n,\dots,M_1} : \bigoplus_p \ (M_n \otimes\cdots \otimes M_1)^G_{p, \mathsf{irr}}\longrightarrow (M_n\otimes \cdots \otimes M_1)^G.$$

Next we will prove that any invariant is a sum of tensors of irreducible invariants.
\Blm
$\psi_{M_n,\dots,M_1}$ is surjective.
\Elm
\Bdm
We begin by providing a natural filtration of  the  invariants. For each non-ordered partition $p=(n_l, \dots, n_1)$ of $n$ we consider similar maps as above,  but not restricted to   irreducible invariants: let $\varphi_p$ from
$$(M_n\otimes \cdots \otimes M_{n_1 + \dots +n_{l-1} + 1})^{G}\otimes\cdots\otimes (M_{n_1+n_2}\otimes\cdots\otimes M_{n_1+1})^G\otimes (M_{n_1}\otimes \cdots \otimes M_{1})^{G}$$
to
$(M_n\otimes \cdots \otimes M_1)^G .$ Consider the image of the sum of the maps $\varphi_p$ along all the $l$-partitions $p$ of $n$; we call it \textit{the space of $l$-composite invariants} and we denote it $\left[(M_n\otimes \cdots \otimes M_1)^G\right] ^l$. Note that $2$-composite invariants are  the space of composite invariants that we have defined before.
The following holds:
$$0\subseteq\left[(M_n\otimes \cdots \otimes M_1)^G\right] ^n \subseteq\cdots\subseteq \left[(M_n\otimes \cdots \otimes M_1)^G\right]^{l} \subseteq \cdots $$
$$\subseteq \left[(M_n\otimes \cdots \otimes M_1)^G\right] ^2\subseteq \left[(M_n\otimes \cdots \otimes M_1)^G\right] ^1 = (M_n\otimes \cdots \otimes M_1)^G.$$
\normalsize
Observe that the first stage of this filtration verifies
$$\left[(M_n\otimes \cdots \otimes M_1)^{G}\right]^n = M_n^G\otimes \cdots \otimes M_1^G.$$
 Moreover this sub-space of the invariants is already  in the image of $\psi$ since $M_i^G =
(M_i)^G_{\mathsf{irr}}$ as mentioned above.

Assume that $\left[(M_n\otimes \cdots \otimes M_1)^{G}\right]^l$ is contained in $\im\left(\psi_{M_n,\dots,M_1}\right)$ in order to prove that $$\left[(M_n\otimes \cdots \otimes M_1)^{G}\right]^{l-1}\subset \im\left(\psi_{M_n,\dots,M_1}\right).$$  Let $m\in \left[(M_n\otimes \cdots \otimes M_1)^{G}\right]^{l-1}$ and suppose first that $m$ is obtained from a fixed $l-1$ partition.  Each invariant tensor in the sub--strings determined by the partition decomposes as a composite plus an irreducible. This shows that $m$ is the sum of terms of two kinds:
\begin{itemize}
 \item[-]tensors of irreducible invariants, which belong by definition to $\mathsf{Im}(\psi)$.
\item[-]$(l-1)$-tensors which contain at least one composite invariant, so belonging to $\left[(M_n\otimes \cdots \otimes M_1)^{G}\right]^l$ which we have assumed is contained in $\im(\psi).$
\end{itemize}
A general $m$ is a sum of terms as above.\qed
\Edm

\Bpo Let $M$ be a
  $kG$-module which is finite dimensional as a vector space, let $n$ be a positive  integer and consider the constant sequence $M_n=\cdots =M_1= M$. Then  $\psi_{M,\dots, M}$ is bijective.
\Epo
\Bdm
Let $T_k(M)^G = k\oplus M^G\oplus (M\otimes M)^G \oplus \cdots $ be  the invariant subalgebra of the tensor algebra $T_k(M)$. The theorem of Kharchenko-Lane states that $T_k(M)^G \cong T_k(U)$ where $U$ is a homogeneous sub-vector space of $T_k(M)^G$.

We assert that $U_n = U\cap (M^{^{\otimes n}})^G$ is a vector space of irreducible invariants for every $n$, namely  it is a complement of the canonical subspace of composite invariants. This is clear at the first degree where the invariants are non-zero, since for this degree the space of composites is zero. Since $T_k(U)$ is a tensor algebra, the elements of the space of composite invariants of degree $n$ are the sums of tensors of homogeneous elements of $U$. We assume that $U_i$ is a space of irreducible invariants in degrees less than $n$.  According to the previous result, in degree $n$ the composites are sums of tensors of irreducible of lower degree, that is of elements of $U_i$ for $i<n$, nothing else is reached through composites. Moreover the intersection of the composites with $U_n$ is zero since $T_k(U)$ is free. We infer that $U_n$ is a complement of the space of composite invariants.

The Theorem of Kharchenko-Lane states that the map $T_k(U)\to T_k(M)^G$ is an isomorphism, then $\psi_{M,\dots,M}$ is injective as well. \qed

\Edm

\begin{theo}
Let $M_n,\dots, M_1$ be a sequence of $kG$-modules which are finite dimensional as vector spaces. Then $\psi_{M_n,\dots,M_1}$ is bijective.
\end{theo}
\begin{proo}
Let $M=M_1\oplus\cdots\oplus M_n$. Then the direct sum of the maps  $\psi_{M_{i_n},\dots,M_{i_1}}$ along all the sequences $(i_n,\dots, i_1)$ of
integers belonging to $\{1,\dots,n\}$ equals $\psi_{M,\dots,M}$ which is a bijection. Hence all those maps are invertible, in particular the one
corresponding to the sequences $(n,\dots,1)$.\qed
\end{proo}

\section{\sf Invariants of a free linear category}

Let $\L_k(V)$ be a free $k$-category on a set of objects $\L_0$ with respect to a family $V$ of vector spaces $\{{}_yV_x\}_{x,y\in\L_0}$.

Let $G$ be a finite group acting linearly on each ${}_yV_x$, in other words $V$ is a family of $kG$-modules. Using the universal property of $\L_k(V)$ we infer an action of $G$ by invertible endofunctors of $\L_k(V)$ which are the identity on objects. Note that the resulting action on morphisms is diagonal on tensor products of vector spaces of the family, in other words the action on morphisms is given by the $kG$-structure of a tensor product of $kG$-modules. Our purpose in this Section is to prove that the invariant category $\L_k(V)^G$ is a free $k$-category.

 Let $Q$ be the track-quiver of $\L_k(V)$. Consider $\gamma$ an oriented path in $Q$ of length $n$, that is a sequence of $n$ concatenated arrows
$u_n\leftarrow u_{n-1}\leftarrow\cdots\leftarrow u_0$. Let
$$V_\gamma= {}_{u_n}V_{u_{n-1}}\otimes\cdots \otimes {}_{u_1}V_{u_{0}}$$
be the vector space corresponding to $\gamma$. {We denote ${}_yQ_x$ the set of oriented paths of $Q$ from $x$ to $y$. Then
$$ \li{y}{\left[\L_k(V)\right]_x}= \bigoplus_{\gamma\in {_{y}\!Q_x}} V_\gamma.$$

Observe that each $V_\gamma$ is a tensor product of a sequence of $kG$-modules as considered in the previous Section.  The invariant sub-category $\L_k(V)^G$ has set of objects $\L_0$, and morphisms as follows:
$$\li{y}{\left[\L_k(V)^G\right]_x}= \bigoplus_{\gamma\in {_{y}\!Q_x}} V_\gamma^G.$$
Let $\left[V_\gamma^G\right]^2$ be the space of composite invariants as defined before, and let $\left[V_\gamma^G\right] _{\mathsf{irr}}$ be a vector space of irreducible invariants, that is a chosen subvector space such that $$V_\gamma = \left[V_\gamma^G\right] _{\mathsf{irr}}\ \oplus\ \left[V_\gamma ^G\right]^2.$$

Let $U=\{{}_yU_x\}_{x,y\in\L_0}$ be the family of vector spaces given by }
$${}_yU_x = \bigoplus_{\gamma\in {}_y\!Q_x} \left[V_\gamma^G\right]_{\mathsf{irr}}.$$

\begin{theo}
Let $G$ be a finite group and let $V=\{_yV_x\}_{x,y\in \L_0}$ be a family of $kG$-modules which are finite dimensional as vector spaces.
The functor  $F:\L_k(U) \longrightarrow \L_k(V)^G$ given by the universal property of $\L_k(U)$ is an isomorphism of categories.
 \end{theo}

 \begin{proo}
 Let $R$ be the track-quiver of the free $k$-category $\L_k(U).$ The arrows of $R$ are precisely the oriented paths $\gamma$ of $Q$ such that $\left[V_\gamma^G\right]_{\mathsf{irr}}\neq 0$. Let $x$ and $y$ be objects. Our purpose is to prove that the map
 $${}_yF_x:{}_y\left[\L_k(U)\right]_x\longrightarrow \li{y}{\left[\L_k(V)^G\right]}_x$$
 is bijective.

  Recall that $\L_k(U)$ is free, hence ${}_y\left[\L_k(U)\right]_x$ is the direct sum of the vector spaces $\L_k(U)_\beta$ along the oriented paths $\beta= \beta_l\dots\beta_1\in {}_yR_x$  where the $\beta_i$'s are arrows in $R$, that is the $\beta_i$'s are oriented paths in $Q$ such that $\left[V^G_{\beta_i}\right]_{\mathsf{irr}}\neq 0$.
 Then
 $$\L_k(U)_\beta =\left[V^G_{\beta_l}\right]_{\mathsf{irr}}\otimes\cdots\otimes  \left[V^G_{\beta_1}\right]_{\mathsf{irr}}.$$
 For a given oriented path $\gamma\in{}_yQ_x$ of length $n$ we consider its \textit{partitions} which are in one-to-one correspondence with the non-ordered partitions  $(n_l,\cdots,n_1)$ of $n$ considered in the previous Section. More precisely a partition of $\gamma$ by non trivial sub-paths is $(\beta_l,\dots,\beta_1)$ where $\beta_i$ is a sub-path of $\gamma$ of length  $n_i$.   If
$\left[V^G_{\beta_i}\right]_{\mathsf{irr}}\neq 0$ for all $i$, then there is an oriented path $\beta\in{}_yR_x$  where the $\beta_i$'s are arrows of $R$. For a given $\gamma$, let $B_\gamma$ be the set of $\beta$'s obtained this way. We put $$\L_k(U)_\gamma =\bigoplus_{\beta\in\B_\gamma} \L_k(U)_\beta.$$
The map ${}_yF_x$ decomposes as a direct sum of maps ${}_yF_x(\gamma)$ along $\gamma\in{}_yQ_x$, where
$${}_yF_x(\gamma): \L_k(U)_\gamma\longrightarrow V_\gamma^G.$$
The last Theorem of the previous Section states that each of those maps is bijective.

 \qed
 \end{proo}

\section{\sf Invariants and representation type}
Let $\C$ be a $k$-category. We consider the abelian $k$-category of $\C$-modules, which objects are $k$-functors from $\C$ to the category of finite dimensional vector spaces. Note that for an object-finite $k$-category $\C$, usual left $\oplus\C$-modules of finite dimension coincides with  $\C$-modules.  Observe that the $k$-algebra $k(t)$ has no finite dimensional modules apart from $0$, hence the same holds for the $k$-category with a single object and endomorphisms  $k(t)$.

A  $\C$-module $\M$ is called  \textit{indecomposable} if the only decomposition of $\M$ as a direct sum of two submodules is trivial, \textit{i.e.} $\M\oplus 0$.{The module $\M$ is \textit{simple} if it has no submodules besides the trivial ones, $0$ and $\M$.}

We recall that a \emph{skeleton} of a $k$-category $\C$ is any full sub-category determined by the choice of an object in each isomorphism class of objects of $\C$. A category  and its skeletons are equivalent, consequently the categories of modules are also equivalent. A $k$-category $\C$ is called \emph{skeletal} if it coincides with its skeleton, in other words different objects of $\C$ are not isomorphic. Observe that if $\C$ is skeletal, any full subcategory of $\C$ is also skeletal.

It is straightforward to verify that free $k$-categories are skeletal.

\begin{defi}
A  $k$-category $\C$  is of \emph{finite representation type} if the number of isomorphism classes of indecomposable $\C$-modules is finite.
\end{defi}

Recall that a \emph{two-sided ideal}  $I$ of a a $k$-category $\C$ is a family of sub-vector spaces $_y I_x\subset {_y\C_x}$ such that for every triple of objects $x,y,z$ it  holds that $_z I_y\  {_y\C_x} \subset {_z\C_x}$ and $_z \C_y\  {_yI_x} \subset {_z\C_x}$. The \emph{quotient category} $\C/I$ has the same objects than $\C$ and the morphisms are the quotient vector spaces ${_y\C_x}/{_y I_x}$. Note that the composition is well defined precisely because $I$ is a two-sided ideal.

\begin{exam}
Let $\L= \L_k(V)$ be a free $k$-category over a family of vector spaces $V=\{_y\!V_x\}_{y, x \in \L_0}$. Let  $\left({{}_y{S}_x}\right)^{\geq i}$ be the set of sequences $(y=u_n, u_{n-1}, \cdots, u_0= x)$ where $n\geq i$, and let $W^i$ be the family
 $$ {\li{y}{W}_x}^i = {\bigoplus_{u \in \left({\li{y}{S}_x}\right)^{\geq i}}} W(u).$$
 Then $W^i$ is a two-sided ideal of $\L$ for every $i$. Moreover $W^i = \left(W^1\right)^i$.
\end{exam}

\begin{defi}
A \emph{well presented} $k$-category is a $k$-category of the form $\L/I$ where $\L$ is a free $k$-category, and $I$ is a two-sided ideal of $\L$ contained in $W^1$.
\end{defi}

Observe that a well presented $k$-category is skeletal.

 \begin{lemm} Let $\L/I$ be a well presented $k$-category of finite representation type. Then $\L/I$ is object-finite.
\end{lemm}
\begin{proo}
Let $x$ be a fixed object, and consider the functor $S_x$ from $\L$ to finite dimensional spaces which is given by the universal property (see Proposition \ref{up}) as follows: $S_x(x)=k$ and  $S_x(y) =0$ for $y\neq x$, while the linear maps on the $_yV_x$ are all zero. We observe that the functor $S_x$ is zero on the two-sided ideal $W^1$, hence on $I$. Consequently $S_x$ is a $\L/I$-module. Those modules are clearly simple and they are non isomorphic for different objects. Since the category is of finite representation type the number of objects is finite.\qed
\end{proo}

\normalsize

Let $\C$ be a $k$-category such that that its skeleton $ \mathcal S$ is object-finite.   The \emph{dimension} of  a $\C$-module $\M$ is
$\mathsf{dim}_k\M = \sum_{x\in\mathcal S_0}\mathsf{dim}_k\M(x)$. Note that this integer do not depend on the choice of the skeleton.

\begin{defi}
Let $\C$ be a finite-objects $k$-category.   We say that $\C$ is  of
\emph{tame representation type}  if for each positive integer $d$  the isomorphism classes of indecomposables $\mathcal C$-modules of $k$-dimension bounded by $d$ are contained in a finite number of one-parameter families, except possibly a finite number of them.
\end{defi}

 We recall Gabriel's Theorem, see \cite{ga1972,begepo}: a  free $k$-category $\L_k(V)$ is of finite representation type if and only if the underlying graph of a bases-quiver is a finite disjoint union of Dynkin diagrams $A_n (n\geq 1), D_n (n\geq 4), E_6, E_7$ or $E_8$. This result has been generalized by Donovan and  Freislich in \cite{dofr} as follows: $\L_k(V)$ is of tame representation type if and only if the underlying graph is a  a finite disjoint union of Dynkin or extended Dynkin diagrams $\widetilde{A_n} (n\geq 1), \widetilde{D_n} (n\geq 4), \widetilde{E_6}, \widetilde{E_7}$ or $\widetilde{E_8}$.

\begin{defi}
 Let $\C$ be a $k$-category. $\C$ is\emph{ finitely-objects of finite, respectively tame, representation type} if any object-finite full  $k$-subcategory of $\C$  is respectively of finite or tame representation type.
 \end{defi}

\begin{rema}
In case $\C$ is a finite-object $k$-category which is of finite or tame representation, note that the preceding definition is consistent. Indeed, in that case we assert that any full $k$-subcategory $\B$ of $\C$ is respectively of finite or tame representation type. This is so since a $\B$-module can be extended to a $\C$-module by assigning  zero vector spaces to the objects which are not in $\B$. This provides a functor which  sends non isomorphic  $\B$-indecomposable module to non isomorphic $\C$-indecomposable module with the same $k$-dimension.
\end{rema}
 \begin{exam}
   The free $k$-category which bases-quiver is any orientation of the graph $$\cdots\ \--\bullet\!\--\!\bullet\!\--\!\bullet\!\--\!\bullet\!\--\!\bullet\--\ \cdots$$ is finitely-objects of finite representation type.
   \end{exam}

Note that in general a $k$-subcategory of a free $k$-category of finite representation type is not of tame representation type, as the following well known example shows (compare with Example \ref{notfreesubcat}).
\begin{exam}
Let $\L$ be the free $k$-category whose bases-quiver is
$$u_0\longrightarrow u_1 \longrightarrow u_2 \longrightarrow u_3\longrightarrow u_4\longrightarrow u_5\longrightarrow u_6$$
Consider the subcategory $\B$ with the same objects while
$$0={}_{u_2}\B_{u_1} ={}_{u_3}\B_{u_1}={}_{u_4}\B_{u_1}={}_{u_5}\B_{u_1}={}_{u_3}\B_{u_2}={}_{u_4}\B_{u_2}=$$$${}_{u_5}\B_{u_2}={}_{u_4}\B_{u_3}={}_{u_5}\B_{u_3}={}_{u_5}\B_{u_4}$$
and ${}_y\B_x = {}_y\L_x$ otherwise. Then
$\B$ is the incidence $k$-category of the partially ordered set $\{u_0,u_1,u_2,u_3, u_4,u_5,u_6\}$ where $u_0$ is the smallest element, $u_6$ the largest one, and $u_1,u_2,u_3,u_4$ and $u_5$ are not comparable. This category is not of tame representation since its quotient by the square of the Jacobson radical is not of tame representation type. Indeed its separate quiver is is the union of two stars with $5$ branches which is not a disjoint union of extended Dynkin and/or Dynkin diagrams.

\end{exam}

Next we recall results related to cleaving procedures introduced by D.G. Higman \cite{hi} and J. Jans \cite{ja}, see also \cite{bagarosa,ge}.
\begin{defi}
Let $\C$ be a $k$-category and let $\B$ be a $k$-subcategory with the same set of objects. $\B$ is a  cleaving $k$-subcategory of $\C$ if there exists a family of  vector spaces $X=\{_y\X_x\}_{x,y\in \B_0}$ verifying the following:

- for every pair $x,y$ of objects $_y\X_x\oplus  {} _y\B_x = {}_y\C_x$ ,

- the family provides a $\B$-bimodule, \textit{i.e.} for any triple of objects $x,y,z$
$$_z\B_y\  {}_y\X_x \subseteq  {}_z\X_x \mbox{ and } _z\X_y\  {}_y\B_x \subseteq {}_z\X_x.$$
Then $\X$ is called a cleaving bimodule of $\B$ in $\C$. If the vector spaces of the family $\X$ are finite dimensional, we say that $\B$ is a cofinite cleaving subcategory and that $\X$ is a finite cleaving bimodule of $\B$ in $\C$.
\end{defi}

Let $\M$ be an indecomposable $\C$-module, let $[\M]$ be its isomorphism class and let $\mathsf{Ind}_\B[\M]$ be the set of isomorphism classes of indecomposables modules which are direct summands of the restriction of $\M$ to $\B$.

\begin{lemm}
  Let $\C$ be a finite-object $k$-category and let $\B$ be a cofinite cleaving $k$-sub\-ca\-te\-gory of $\C$. Let $\N$ be an indecomposable $\B$-module.

  There exists an indecomposable $\C$-module $\M$ such that $[\N] \in \mathsf{Ind}_\B[\M]$. Moreover there exists a positive integer $\gamma$ dependent only on $\B$ and $\C$ such that $\M$ can be chosen to verify $$\mathsf{dim}_k\M\leq\gamma\mathsf{dim}_k\N.$$

\end{lemm}
\begin{proo}
Since the categories involved are object-finite, it is equivalent to prove the result for a $k$-algebra $C$ containing a $k$-subalgebra $B$ and a finite dimensional $B$-bimodule $X$ such that $C=B\oplus X$.

 Let $N$ be an indecomposable $B$-module which is finitely dimensional over $k$,  and let  $C\otimes_{B}N$ be the extended $C$-module.  Restricting  to $B$ provides a decomposition into a direct sum of $B$-modules
 $$C\otimes_{B}N = N \ \oplus \ \left(X\otimes_B N\right).$$
  The preceding equality shows that all the modules are finite dimensional over $k$, hence indecomposable modules have local endomorphisms rings and the Krull-Remak-Schmidt theorem is in force. We infer that the isomorphism class of $N$ appears in any decomposition of $C\otimes_{B}N$ into indecomposable $B$-modules. Consider a decomposition $C\otimes_{B}N= \oplus M_i$ into  indecomposable $C$-modules and next a decomposition of the restriction of each $M_i$ to $B$ into a direct sum of indecomposable $B$-modules. We infer that there exists at least an $M_i$ such that $N\in\mathsf{Ind}_B [M_i]$.

 Moreover the $k$-dimension of $M$ is at most $\mathsf{dim}_k N + \mathsf{dim}_k X\mathsf{dim}_k N$, which shows that $\gamma$ can be chosen to be equal to $1+\mathsf{dim}_k X.$
\qed
\end{proo}

    \begin{theo}\label{reptypecofcleaving}
Let $\C$ be a $k$-category and let $\B$ be a cleaving cofinite $k$-sub\-ca\-te\-go\-ry of $\C$. If $\C$ is  finitely-objects of finite or tame representation type, then $\B$ is finitely-objects of respectively finite or tame representation type.
\end{theo}

\begin{proo} We first show that the result follows if it holds for finite-objects $k$-categories. Let $\C$ be a $k$-category finitely-objects of finite or tame representation type, and let $\X$ be a finite cleaving bimodule of $\B$ in $\C$. Let $\B'$ be a full subcategory of $\B$ determined with a finite set of objects $\B'_0$. Let $\C'$ be the full subcategory of $\C$ with set of objects $\B'_0$ and let $\X'$ be the subfamily of $\X$ given by the pairs of objects in $\B'_0$. Note that $\X'$  is a finite cleaving bimodule of $\B'$ in $\C'$. Since we assume that the result holds for finite-object categories, we infer the required conclusion for $\B'$.

We assume now that $\C$ has a finite number of objects and $\B$ is a cleaving cofinite $k$-subcategory.  Let $U$ be the union of the sets  $\mathsf{Ind}_\B[\M]$ where $[\M]$ varies among the isomorphism classes of  indecomposable $C$-modules. The Lemma above shows that $U$ is the set of isomorphism classes of indecomposable $\B$-modules. If $\C$ is of finite representation type, $U$ is finite, hence $\B$ is of finite representation type.

Let $\mathsf{Ind}_d\B$ be the set of isomorphism classes of indecomposable $\B$-modules of dimension at most $d$. Using the Lemma above, we observe
that $\mathsf{Ind}_d\B$  is contained in the subset of $U$ of isomorphism classes of indecomposable $\C$-modules of dimension at most $\gamma d$. If
$\C$ is of tame representation type, this subset is parametrized  by a finite number of copies of $k$ apart a finite number of isomorphism classes.
Consequently this is also the case for $\mathsf{Ind}_d\B$.

\qed
\end{proo}
Recall that a  $k$-category is said to be \textit{hom-finite} if its vector spaces of morphisms are finite dimensional. Next we prove a Theorem by I. Reiten and C. Riedtmann in \cite{reri}, which is valid for an arbitrary hom-finite $k$-category with an action by a group $G$ which order is invertible in the field $k$. This way we obtain a proof for the Theorem \ref{Thethm} for the particular case where the free $k$-category is hom-finite and  $|G|\neq 0$ in $k$.

\begin{theo}
Let $G$ be a finite group which order is not zero in the field $k$. Let $\C$ be a hom-finite $k$-category and let $G$ be a group acting on $\C$ with trivial action on the objects.
If $\C$ is  finitely-objects of finite or tame representation type,  $\C^G$ is finitely-objects of respectively finite or tame representation type.
\end{theo}

\begin{proo}
We assert that $\C^G$ is a cofinite cleaving subcategory of $C$. Indeed each space of morphisms of $\C$ is a finite dimensional $kG$-module. The hypothesis on the order of $G$ insures that the group algebra $kG$ is semisimple. Then $${}_y\C_x= [{}_y\C_x]^G \oplus {}_y\X_x$$ where ${}_y\X_x$ is the canonical complement, namely is the direct sum of the isotypic components of non trivial simple $kG$-submodules. Clearly this family constitutes a $\C^G$-bimodule.
\qed
\end{proo}

We will need the following result which is certainly well known, we provide a proof for completeness.

\begin{lemm}\label{fullfree}
Let $\L_k(V)$ be a free $k$-category and let $\C$ be a full subcategory of it. Then $\C$ is also a free $k$-category.
\end{lemm}
\begin{proo}
We first prove the analogous result for small categories, i.e. without $k$-structure. Let $Q$ be a quiver and let $\F_Q$ be the associated free category (see Remark \ref{free}). Let $\C$ be a full subcategory of $\F_Q$. A morphism of $\C$ is a sequence of concatenated arrows of $Q$ (i.e. a path) having source and target objects  in $C_0$. A path of $\C$ is called $\C$-\textit{primitive} if the vertices of all of its arrows are not in $\C_0$ - except of course the source and the target of the path. Let $R$ be the quiver with vertices $\C_0$ and  arrows $a_w$ associated to all $\C$-primitive paths $w$,  the source and target objects of $a_w$ are the same than for $w$. There is an evident functor from $\F_R$ to $C$ determined by the universal property of free categories, which is the identity on objects and sends each arrow $a_w$ of $R$ to the  $\C$-primitive path $w$. The functor is bijective on morphisms since any morphism of $\C$ can uniquely be decomposed as a concatenation of $\C$-primitive morphisms, just by cutting the morphism at the vertices of the arrows which are in $\C_0.$ Consequently $C$ is isomorphic to the free $k$-category $\F_R$.

Let now $\C$ be a full $k$-subcategory of a free $k$-category $\L_k(V)$. Let $Q$ be the bases-quiver of the latter. Let $\F_Q$ be the free category and let $\widehat{\C}$ be the full subcategory of $\F_Q$ determined by $\C_0$. According to the above consideration $\widehat{\C}$ is a free category. We observe that its linearization $k\widehat{\C}$ is the full $k$-subcategory of $\C$ with objects $\C_0$. We infer $k\widehat{\C}=\C$. Then $\C$ is the linearization of a free category, hence is a free $k$-category.
\qed
\end{proo}

\begin{prop}\label{sgcleaving}
Let $\L_k(V)$ be a Schurian-generated free $k$-category (see De\-fi\-nition \ref{sg}). Let $G$ be a finite group acting homogeneously on it. Then $\L_k(V)^G$ is a cleaving subcategory of $\L_k(V)$.
\end{prop}

\begin{proo}
As in the proof of Proposition \ref{sgfree}, let $\chi_w$ be the character associated to each path $w$ of positive length of the bases-quiver. $\L_k(V)^G$ is the linearization of the category given by the paths $w$ such that $\chi_w=1$. Let $\X$ be the family of sub-vector spaces of the morphisms with bases the paths with non-trivial characters. Clearly $\X$ provides a $\L_k(V)^G$-bimodule. Moreover all the paths from $x$ to $y$ form a basis of ${}_y[\L_k(V)]_x$, hence we have a direct sum decomposition. This shows that $\X$ is a cleaving bimodule.
\qed
\end{proo}

\begin{theo}\label{Thethm}
Let $\L$ be a free $k$-category over a set $\L_0$ and over a family of $k$-vector spaces $V=\{{}_yV_x\}_{x,y\in\L_0}$. Let $G$ be a finite group acting homogeneously on $\L$, equivalently each vector space of the family is provided with a structure of a $kG$-module. If $\L$ is   finitely-objects of finite or tame representation type, $\L^G$ is finitely-objects of respectively finite or tame representation type.
\end{theo}
\begin{proo}First we show that it is enough to prove the result for a finite-object $k$-category. Indeed let $\K$ be a finite-object full subcategory of $\L^G$. Let $\K'$ be the full subcategory of $\L$ determined by the finite set $\K_0$. We recall that the action of $G$ is trivial on objects, hence $\K'$ has also an action of $G$ and $\K'^G=\K$. Note that $\K'$ is free according to the Lemma \ref{fullfree}. Since its number of objects is finite, by assumption $\K'$ is respectively of finite or tame representation type. If the theorem is proved for finite-object categories we infer that $\K'^G$, that is $\K$, is of finite respectively tame representation type.

We assume now that $\L$ is a finite-object free $k$-category, so its bases-quiver is a union of Dynkin or extended Dynkin graphs with some orientation. We will analyze separately  the two following cases:
\begin{itemize}
\item[-]$\widetilde A_n$  with cyclic orientation, that is an oriented crown as in Example \ref{crown}. Note that the corresponding free $k$-category is not hom-finite,
\item[-]$\widetilde A_1$ with non cyclic orientation, that is the Kronecker quiver $\cdot\rightrightarrows\cdot$ Note that the corresponding free  $k$-category is not Schurian-generated.
\end{itemize}
In all other cases the corresponding free $k$-categories are Schurian-generated and hom-finite, consequently the invariant sub-category is cleaving cofinite according to the previous Proposition. The result follows from Theorem \ref{reptypecofcleaving}.

Let $Q_n$  be an oriented crown with $n$ vertices, that is the graph $\widetilde A_{n-1}$ with a cyclic orientation  provided with an action of a finite group $G$ which is trivial on vertices. In other words, to each arrow there is an attached character $G\to k^{\cdot}$. Our purpose is to prove that the invariant $k$-category is of tame representation type. Note that the free $k$-category determined by $Q_n$ is Schurian-generated, then as in the proof of Proposition \ref{sgfree} each path $w$ has an associated character. Recall from the proof of this Proposition that the irreducible invariants paths are the arrows of the bases-quiver of the invariant category. Given a vertex of $Q_n$, we assert that there is precisely one irreducible invariant path having this object as a source. Indeed, the group of characters is finite abelian, let $e$ be its exponent. Any product of $e$ characters is trivial, hence there is at least one invariant path with a given source. Then there is at least one irreducible invariant path with given source. Finally observe that for an oriented crown, two paths $w$ and $w'$ with a common source are either equal or there is a path $w''$ such that $w'=w''w$ or $w=w''w'$. Hence if $w$ and $w'$ are invariants paths then $w''$ is also an invariant path. Consequently there is a unique irreducible invariant path with a given source. Similarly we assert that there is precisely one irreducible path with a given target.

The bases-quiver of the invariant category inherits this property, namely each vertex is the source of a unique arrow and the target of a unique arrow. It follows that the quiver is a union of oriented crowns.

Finally consider the Kronecker quiver. An action of a finite group $G$ which is trivial on objects is given by a $kG$-module $V$ of dimension $2$. According to the dimension  of the submodule $V^G$ of $V$, equal to $2$, $1$ or $0$, the respective invariants are given by $\widetilde A_1$, $A_2$ or two vertices.
\qed
\end{proo}

\footnotesize
\noindent C.C.:
\\Institut Montpelli\'{e}rain Alexander Grothendieck (IMAG) UMR 5149\\
Universit\'{e}  de Montpellier, F-34095 Montpellier cedex 5,
France.\\
{\tt Claude.Cibils@umontpellier.fr}

\medskip

\noindent E.N.M.:\\
Departamento de Matem\'atica, IME-USP,\\
Rua do Mat\~ao 1010, cx pt 20570, S\~ao Paulo, Brasil.\\
{\tt enmarcos@ime.usp.br}

\end{document}